\documentclass[12pt]{article}
\usepackage{amsmath,amssymb,amsthm}
\usepackage{cite}

\def\barr{\begin{array}}
\def\earr{\end{array}}

\title{A new Menon-type identity derived from group actions}
\author{Marius T\u arn\u auceanu}
\date{May 26, 2021}

\begin{document}

\maketitle

\begin{abstract}
In this short note, we give a new Menon-type identity involving the sum of element orders and the sum of cyclic subgroup orders of a finite group.
It is based on applying the weighted form of Burnside's lemma to a natural group action.
\end{abstract}
\smallskip

{\small
\noindent
{\bf MSC 2020\,:} Primary 11A25; Secondary 20D60.

\noindent
{\bf Key words\,:} Menon's identity, weighted Burnside's lemma, group action, sum of element orders, sum of cyclic subgroup orders.}

\section{Introduction}

One of the most interesting arithmetical identities is due to P.K. Menon \cite{5}.

\bigskip\noindent{\bf Menon's identity.}
{\it For every positive integer $n$ we have
\begin{equation}
\sum_{a\in\mathbb{Z}^*_n}\gcd(a-1,n)=\varphi(n)\,\tau(n),\nonumber
\end{equation}where $\mathbb{Z}^*_n$ is the group of units of the ring $\mathbb{Z}_n=\mathbb{Z}/n\mathbb{Z}$, $\gcd(,)$ represents
the greatest common divisor, $\varphi$ is the Euler's totient function and $\tau(n)$ is the number of divisors of $n$.}
\bigskip

There are several approaches to Menon’s identity and many generalisations. One of the methods used to prove Menon-type identities is based on the Burnside's Lemma concerning group actions (see e.g. \cite{5,6,7,8,9,10}). In what follows, we will use a generalization of this result, called the Weighted Form of Burnside's Lemma (see e.g. \cite{2}).\newpage

\noindent{\bf Weighted Form of Burnside's Lemma.}
{\it Given a finite group $G$ acting on a finite set $X$, we denote
\begin{equation}
Fix(g)=\{x\in X : g\circ x=x\},\, \forall\, g\in G.\nonumber
\end{equation}Let $R$ be a commutative ring containing the rationals and $w:X\longrightarrow R$ be a weight function that is constant on the distinct orbits $O_{x_1}$, ..., $O_{x_k}$ of $X$. For every $i=1,...,k$, let $w(O_{x_i})=w(x)$, where $x\in O_{x_i}$. Then
\begin{equation}
\sum_{i=1}^k \,w(O_{x_i})=\frac{1}{|G|}\,\sum_{g\in G}\sum_{\,x\in Fix(g)}\!\!w(x).
\end{equation}}

Note that the Burnside's Lemma is obtained from (1) by taking the weight function $w(x)=1$, $\forall\, x\in X$.
\smallskip

Next we will consider a finite group $G$ of order $n$ and the functions
\begin{equation}
\psi(G)=\sum_{g\in G}\,o(g)\, \mbox{ and }\, \sigma(G)=\!\!\!\!\sum_{H\in C(G)}|H|,\nonumber
\end{equation}where $o(g)$ is the order of $g\in G$ and $C(G)$ is the set of cyclic subgroups of $G$.\,\footnote{\,For more details concerning these functions, we refer the reader to \cite{1,3} and \cite{4}, respectively.} Also, for every divisor $m$ of $n$, we will denote $G_m=\{g\in G : g^m=1\}$.\vspace{3mm}

Our main result is stated as follows.

\bigskip\noindent{\bf Theorem 1.}
{\it Under the above notations, we have
\begin{equation}
\sum_{a\in\mathbb{Z}^*_n}\psi(G_{\gcd(a-1,n)})=\varphi(n)\,\sigma(G).\vspace{-2mm}
\end{equation}}
\smallskip

Clearly, (2) gives a new connection between the above functions $\psi(G)$ and $\sigma(G)$. We remark that an alternative way of writing (2) is
\begin{equation}
\sum_{a\in\mathbb{Z}^*_n}\sum_{\,d|\gcd(a-1,n)}d\varphi(d)n_d(G)=\varphi(n)\,\sigma(G),
\end{equation}where $n_d(G)$ denotes the number of cyclic subgroups of order $d$ in $G$, for all $d$ dividing $n$.
\smallskip\newpage

For $G=\mathbb{Z}_n$, Theorem 1 leads to the following corollary.

\bigskip\noindent{\bf Corollary 2.}
{\it We have
\begin{equation}
\sum_{a\in\mathbb{Z}^*_n}\psi(\mathbb{Z}_{\gcd(a-1,n)})=\varphi(n)\,\sigma(n),
\end{equation}where $\sigma(n)$ is the sum of divisors of $n$.}
\bigskip

Finally, since $\psi(\mathbb{Z}_n)\geq\frac{q}{p+1}\,n^2$, where $q$ and $p$ are the smallest and the largest prime divisor of $n\geq 2$ (see the proof of Lemma 2.9(2) in\cite{3}), from (4) we infer the following inequalities.

\bigskip\noindent{\bf Corollary 3.}
{\it We have
\begin{equation}
\frac{q}{p+1}\,\frac{1}{\varphi(n)}\sum_{a\in\mathbb{Z}^*_n}\gcd(a-1,n)^2\leq\sigma(n)\leq\frac{1}{\varphi(n)}\sum_{a\in\mathbb{Z}^*_n}\gcd(a-1,n)^2.
\end{equation}}\vspace{-3mm}

\section{Proof of Theorem 1}

Let $\mathbb{Z}^*_n=\{a\in\mathbb{N} : 1\leq a\leq n, \gcd(a,n)=1\}$ be the group of units (mod $n$). The natural action of $\mathbb{Z}^*_n$ on $G$ is defined by
$$a\circ g=g^a,\, \forall\, (a,g)\in \mathbb{Z}^*_n\times G.$$Then two elements of $G$ belong to the same orbit if and only if they generate the same cyclic subgroup. This shows that the weight function $w:G\longrightarrow\mathbb{R}$, $w(g)=o(g),\, \forall\, g\in G$, is constant on the distinct orbits $O_{g_1}$, ..., $O_{g_k}$ of $G$. Thus we can apply the  Weighted Form of Burnside's Lemma.

First of all, we observe that $w(O_{g_i})=o(g_i)=|\langle g_i\rangle|$, $\forall\, i=1,...,k$, and therefore the left side of (1) is $\sigma(G)$. 

Next we will prove that $Fix(a)=G_{\gcd(a-1,n)}$, for any $a\in \mathbb{Z}^*_n$. Indeed, if $g\in Fix(a)$ then $g^a=g$, that is $g^{a-1}=1$. Since $|G|=n$, we also have $g^n=1$. Consequently, $g^{\gcd(a-1,n)}=1$, i.e. $g\in G_{\gcd(a-1,n)}$. The converse inclusion is obvious.

Now (1) becomes
\begin{equation}
\sigma(G)=\frac{1}{\varphi(n)}\,\sum_{a\in\mathbb{Z}^*_n}\sum_{\,g\in G_{\gcd(a-1,n)}}\!\!o(g)=\frac{1}{\varphi(n)}\,\sum_{a\in \mathbb{Z}^*_n}\psi(G_{\gcd(a-1,n)}),\nonumber
\end{equation}as desired.\hfill\rule{1,5mm}{1,5mm}

\vspace*{3ex}\small

\hfill
\begin{minipage}[t]{5cm}
Marius T\u arn\u auceanu \\
Faculty of  Mathematics \\
``Al.I. Cuza'' University \\
Ia\c si, Romania \\
e-mail: {\tt tarnauc@uaic.ro}
\end{minipage}

\end{document}